\documentclass[runningheads]{llncs}

\usepackage{amsmath}
\usepackage{amssymb}

\usepackage{booktabs} 
\usepackage[ruled,vlined]{algorithm2e}
\usepackage{listings}
\usepackage{breakcites}

\usepackage[british]{babel}

\usepackage{tikz}
\usepackage{pgfplots}
\pgfplotsset{compat=newest}
\usepgfplotslibrary{colormaps}

\DeclareRobustCommand{\stirling}{\genfrac\{\}{0pt}{}}

\usepackage{graphicx}

\usepackage[colorlinks]{hyperref}
\usepackage[nameinlink]{cleveref}

\hypersetup{
  colorlinks=true,
  citecolor=green
}

\begin{document}
\title{Handy Formulas for Binomial Moments}
%
%
\author{Maciej Skorski\inst{1}}
\authorrunning{M. Skorski}
%
\institute{University of Luxembourg\\
\email{maciej.skorski@gmail.com}}
\maketitle              
\begin{abstract}
Despite the relevance of the binomial distribution for probability theory and applied statistical inference, its higher-order moments
are poorly understood. The existing formulas are either not general enough, or not structured and simplified enough for intended applications. 

This paper introduces novel formulas for binomial moments, in form of \emph{polynomials in the variance} rather than in the success probability. The obtained formulas are arguably better structured, simpler and superior in their numerical properties compared to prior works. In addition, the paper presents algorithms to derive these formulas along with working implementation in the Python symbolic algebra package.

The novel approach is a combinatorial argument coupled with clever algebraic simplifications which rely on symmetrization theory. As an interesting byproduct we establish \emph{asymptotically sharp estimates for central binomial moments}, improving upon partial results from prior works.

\keywords{Binomial Distribution  \and High-order Moments \and Symbolic Algebra.}
\end{abstract}
\section{Introduction}

\subsection{Background and Related Work}

The binomial distribution $\mathrm{Binom}(n,p)$, which counts the total number of successes within $n$ independent trials each succeeding with probability $p$, is of historical and fundamental importance for probability theory and applied statistical inference. In particular, it
appears in quantitative variants of central limit theorem~\cite{laplace1810memoire,de1733approximatio,uspensky1937introduction}, and is broadly used in statistical modeling~\cite{little1989testing,agresti1998approximate,fraas1994binomial,daniel2003binomial,yang2011new} also as a building block of more advanced models~\cite{young2008pooling}; this includes performing A/B tests on conversion rates popular in business~\cite{goodson2014most}.

Despite this large body of work on approximate inference, little is known about the \emph{exact higher moments} of the binomial distribution. Except being of natural interest, the demand for such formulas comes from seeking for provable guarantees on probability tails by means of the \emph{moment method}; for example, binomial bounds are needed for analyzing random projections~\cite{jagadeesan2019understanding}.

While the textbooks usually cover only the variance, sometimes also the skewness and kurtosis), 
there have been only few research papers discussing formulas for binomial moments of order $d>4$.
The first recursion formula for binomial moments appeared in \cite{benyi2005recursive} for the special case of $p=\frac{1}{2}$.
The case of general $p$ was handled in~\cite{griffiths2013raw} by means of recursions
utilizing Stirling numbers of the first kind. This was subsequently simplified by means of moment generating functions in~\cite{knoblauch2008closed}, and resulted in a more compact formula involving Stirling numbers of the second kind. Very recently, a recursion-free derivation of raw moments has been presented in~\cite{nguyen2019probabilistic}. The common idea is to see the moments as (more or less explicit) polynomials in $n$ and $p$ and group terms to make the formula more compact.

The discussed approaches still do not offer a satisfactory answer, as the formulas are not handy enough to be directly applicable. The author of the most general formula in~\cite{knoblauch2008closed} didn't manage to obtain non-naive bounds on the binomial moments: the bound $O(nq)^d$ with $q=1-p$ on the $d$-th central moment~\cite{knoblauch2008closed} valid for $p<\frac{1}{2}$ is trivial as the centered binomial random variable is bounded between $-p n$ and $q n$ - no extra formulas are needed; this bound is far from the true behavior $O(nq)^{d/2}$ when $nq\to\infty$ (obtained by the Central Limit Theorem). The main formula in~\cite{knoblauch2008closed} is actually a mixture of positive and negative contributions, which makes its numerical convergence problematic and theoretical analysis very difficult (as seen above). Moreover, all the prior works do not exploit the symmetry and produce overly complicated formulas in terms of $p$; it should be noted~\cite{wiki:binomial} that the simplest expressions on the central moments of small orders appear to be obtained with the variance $\sigma^2=p(1-p)$ as a variable. Lastly, the discussed prior works are rather scarce in their presentation of related works and techniques, in particular they seem to have overlooked that the formulas with the appearance of Stirling numbers follow easier by the established approach of factorial moments~\cite{joarder1997classroom,balakrishnan1998note}.

\subsection{Summary of Contributions}
Addressing the aforementioned issues with approaches in prior works, this paper offers the following novel contributions on computing the binomial moments:
\begin{itemize}
    \item \textbf{link to factorial moment} which simplifies the approach from prior works
    \item \textbf{variance-formula} for equivalent yet simpler expressions in $\sigma^2=p(1-p)$
    \item \textbf{algorithm and implementation} for finding the variance-formula
    \item \textbf{stable formula} as explicit sum with positive terms
    \item \textbf{asymptotically sharp bounds} on binomial moments as an application.
\end{itemize}
In summary, when compared to prior works, these results brings a broader scope of the techniques, as well as lead to arguably more handy formula; another added value is the contributed algorithm and its Python implementation\footnote{For code and examples see \url{http://github.com/maciejskorski/binomial_moments}}.

\subsection{Preliminaries}

\subsubsection{Binomial Distribution}
A random variable $S$ follows the binomial distribution with parameters $n$ and $p$, denoted as  $S\sim \mathrm{Binom}(n,p)$,
 when the probability density function is
\begin{align}\label{eq:density}
    \Pr[S=k] = \binom{n}{k}p^kq^{1-k},\quad q\triangleq 1-p,\quad k=0\ldots n.
\end{align}
\subsubsection{Moments}
Let $d$ be a positive integer. The \emph{raw moment} of order $d$ of a random variable $S$ is defined as $\mathbb{E}[S^d]$, while the \emph{central moment} of order $d$ of $S$ equals $\mathbb{E}[(S-\mathbb{E}[S])^d]$. We also use the \emph{factorial moment} defined as $\mathbb{E}[\mathbb{S}^{\underline{d}}]$
where $x^{\underline{d}}=x(x-1)\cdots (x-(d-1))$ is called the $d$-th \emph{falling power}~\cite{graham1989concrete}.
\subsubsection{Special Numbers}
To state some of our results we need \emph{Stirling numbers of second kind}. The symbol $\stirling{n}{k}$ stands for the number of ways of partitioning an $n$ element set into $k$ non-empty subsets. We also need \emph{multinomial coefficients} defined as $\binom{d}{d_1\ldots d_k} = \frac{d!}{d_1!\cdots d_k!}$ when $\sum_{i=1}^{k}d_i=d$ and $\min_i d_i\geqslant 0$ and 0 otherwise, which extend the binomial coefficients. By the multinomial theorem we have that
$(x_1+\ldots+x_n)^d = \sum_{d_1,\ldots,d_n}\binom{d}{d_1\ldots d_n}x^{d_1}\cdots x_n^{d_n}$, extending the binomial formula.

\subsubsection{Polynomials}
To work out the desired polynomial formulas we need some standard algebraic notation. By $\mathbb{Z}[x_1,\ldots,x_k]$ we denote polynomials with integer coefficients in variables $x_1,\ldots,x_k$. A polynomial is symmetric if after exchanging any two variables its sign doesn't change, and anti-symmetric when the sign gets negated.
The fundamental theorem of symmetric polynomials states that any symmetric polynomial from $\mathbb{Z}[x_1,\ldots,x_k]$ can be written as a polynomial in the elementary symmetric functions $s_j(x_1,\ldots,x_n) = \sum_{1\leqslant i_1\ldots \leqslant i_j} x_{i_1}\cdots x_{i_j}$ for $j=1\ldots k$, with integer coefficients. Moreover, anti-symmetric polynomials can be written as a product of a symmetric polynomial and Vandermonde's determinant $\prod_{1\leqslant i<j \leqslant k}(x_i-x_j)$
 (see for example~\cite{prasolov2004polynomials,zhou2003introduction}).


\subsection{Results}
Below we discuss the contributions in more detail, deferring proofs to the end part of the paper. We denote $S\sim \mathrm{Binom}(n,p)$ and fix a positive integer $d$.

\subsubsection{Raw Binomial Moments and Factorial Moments}

Our first result is derivation of a closed-form formula for raw binomial moments. 
This formula appears already in prior works~\cite{benyi2005recursive,knoblauch2008closed,griffiths2013raw},
however our novelty is in the techniques: as opposed to recursion-based approaches  \cite{benyi2005recursive,knoblauch2008closed,griffiths2013raw} we give two alternative proofs a)
by linking central and factorial moments b) by developing a direct counting argument.
In the context of the prior works the approach (a) broadens the perspective and brings pedagogical value,
and the approach (b) will be reused later in the discussion of central moments.

\begin{theorem}[Formula for Raw Binomial Moments]\label{cor:binom_raw}
Then
\begin{align}\label{eq:raw_binom_moment}
\mathbb{E}[S^d] = \sum_{k=0}^{d}n^{\underline{k}}\stirling{d}{k}p^k.
\end{align}
\end{theorem}
The proofs appear respectively in \Cref{proof:binom_raw:factorial} and \Cref{proof:binom_raw:direct}. Below in \Cref{tab:raw_binom} we list the explicit expressions for the first 10 moments.

\begin{table}[!h]
\resizebox{0.99\linewidth}{!}{
\begin{tabular}{ll}
\toprule
$d$ &                                                                                                                                                                                                                                                                                                      $\mathbb{E}[S]^d,\quad S\sim\mathrm{Binom}(n,p)$ \\
\midrule
0 &                                                                                                                                                                                                                                                                $2 p^{2} {\binom{n}{2}} + p {\binom{n}{1}}$ \\
1 &                                                                                                                                                                                                                                       $6 p^{3} {\binom{n}{3}} + 6 p^{2} {\binom{n}{2}} + p {\binom{n}{1}}$ \\
2 &                                                                                                                                                                                                           $24 p^{4} {\binom{n}{4}} + 36 p^{3} {\binom{n}{3}} + 14 p^{2} {\binom{n}{2}} + p {\binom{n}{1}}$ \\
3 &                                                                                                                                                                              $120 p^{5} {\binom{n}{5}} + 240 p^{4} {\binom{n}{4}} + 150 p^{3} {\binom{n}{3}} + 30 p^{2} {\binom{n}{2}} + p {\binom{n}{1}}$ \\
4 &                                                                                                                                                 $720 p^{6} {\binom{n}{6}} + 1800 p^{5} {\binom{n}{5}} + 1560 p^{4} {\binom{n}{4}} + 540 p^{3} {\binom{n}{3}} + 62 p^{2} {\binom{n}{2}} + p {\binom{n}{1}}$ \\
5 &                                                                                                                $5040 p^{7} {\binom{n}{7}} + 15120 p^{6} {\binom{n}{6}} + 16800 p^{5} {\binom{n}{5}} + 8400 p^{4} {\binom{n}{4}} + 1806 p^{3} {\binom{n}{3}} + 126 p^{2} {\binom{n}{2}} + p {\binom{n}{1}}$ \\
6 &                                                                              $40320 p^{8} {\binom{n}{8}} + 141120 p^{7} {\binom{n}{7}} + 191520 p^{6} {\binom{n}{6}} + 126000 p^{5} {\binom{n}{5}} + 40824 p^{4} {\binom{n}{4}} + 5796 p^{3} {\binom{n}{3}} + 254 p^{2} {\binom{n}{2}} + p {\binom{n}{1}}$ \\
7 &                                          $362880 p^{9} {\binom{n}{9}} + 1451520 p^{8} {\binom{n}{8}} + 2328480 p^{7} {\binom{n}{7}} + 1905120 p^{6} {\binom{n}{6}} + 834120 p^{5} {\binom{n}{5}} + 186480 p^{4} {\binom{n}{4}} + 18150 p^{3} {\binom{n}{3}} + 510 p^{2} {\binom{n}{2}} + p {\binom{n}{1}}$ \\
8 &  $3628800 p^{10} {\binom{n}{10}} + 16329600 p^{9} {\binom{n}{9}} + 30240000 p^{8} {\binom{n}{8}} + 29635200 p^{7} {\binom{n}{7}} + 16435440 p^{6} {\binom{n}{6}} + 5103000 p^{5} {\binom{n}{5}} + 818520 p^{4} {\binom{n}{4}} + 55980 p^{3} {\binom{n}{3}} + 1022 p^{2} {\binom{n}{2}} + p {\binom{n}{1}}$ \\
\bottomrule
\end{tabular}
}
\caption{Formulas for Raw Binomial Moments.}
\label{tab:raw_binom}
\end{table}

\subsubsection{Central Binomial Moments}

\paragraph{Symmetric Structure} While in prior works the formulas are derived in terms of $p$, we go beyond that exploiting the symmetry
and showing that the formulas can be written in terms of the variance $\sigma^2=p(1-p)$, which makes them much simpler.
The following theorem proves what can be conjectured by inspection of known formulas for small-order moments~\cite{wiki:binomial}.

\begin{theorem}[Variance-Based Formula]\label{thm:structure}
For $S\sim\mathrm{Binom}(n,p)$ the moment $\mathbb{E}[(S-\mathbb{E}[S])^d]$ is a symmetric polynomial in $p$ and $q$ when $d$ is even, 
and anti-symmetric when $d$ is odd.  In particular denotting $\sigma^2\triangleq pq$ we have
\begin{align}
\mathbb{E}[(S-\mathbb{E}[S])^d] \in
\begin{cases}
  \mathbb{Z}[n,\sigma^2] & d\text{ even} \\
  (1-2p)\mathbb{Z}[n,\sigma^2] & d\text{ odd} 
\end{cases}
\end{align}
\end{theorem}
\Cref{tab:moments} illustrates this result, providing explicit moments of order $d=2\ldots 10$.
The practical usefulness of the formula guaranteed by \Cref{thm:structure} is its simplicity when compared to representation in terms of $p$ alone. The result is intuitive, but not straightforward to prove; we give two arguments based on a) theory of symmetric functions, see
\Cref{proof:thm:structure:algebra}
 and b) our novel combinatorial formula, see \Cref{proof:thm:structure:combinatorics}.
The algorithm deriving the exact formulas is discussed later.

\begin{table}[!ht]
\resizebox{1.0\textwidth}{!}{
\centering
\begin{tabular}{ll}
\toprule
$d$ &                                                                                                                                                                                                                                                                                                                                                                       $\mathbb{E}[(S-\mathbb{E}[S])^d],\quad S\sim\mathrm{Binom}(n,p)$ \\
\midrule
2 &                                                                                                                                                                                                                                                                                                                                                              $n \sigma^{2}$ \\
3 &                                                                                                                                                                                                                                                                                                                                       $n \sigma^{2} \left(- 2 p + 1\right)$ \\
4 &                                                                                                                                                                                                                                                                                                           $3 n^{2} \sigma^{4} + n \left(- 6 \sigma^{4} + \sigma^{2}\right)$ \\
5 &                                                                                                                                                                                                                                                                     $\left(- 2 p + 1\right) \left(10 n^{2} \sigma^{4} + n \left(- 12 \sigma^{4} + \sigma^{2}\right)\right)$ \\
6 &                                                                                                                                                                                                                                    $15 n^{3} \sigma^{6} + n^{2} \left(- 130 \sigma^{6} + 25 \sigma^{4}\right) + n \left(120 \sigma^{6} - 30 \sigma^{4} + \sigma^{2}\right)$ \\
7 &                                                                                                                                                                                               $\left(- 2 p + 1\right) \left(105 n^{3} \sigma^{6} + n^{2} \left(- 462 \sigma^{6} + 56 \sigma^{4}\right) + n \left(360 \sigma^{6} - 60 \sigma^{4} + \sigma^{2}\right)\right)$ \\
8 &                                                                                                                                   $105 n^{4} \sigma^{8} + n^{3} \left(- 2380 \sigma^{8} + 490 \sigma^{6}\right) + n^{2} \left(7308 \sigma^{8} - 2156 \sigma^{6} + 119 \sigma^{4}\right) + n \left(- 5040 \sigma^{8} + 1680 \sigma^{6} - 126 \sigma^{4} + \sigma^{2}\right)$ \\
9 &                                                                                          $\left(- 2 p + 1\right) \left(1260 n^{4} \sigma^{8} + n^{3} \left(- 13216 \sigma^{8} + 1918 \sigma^{6}\right) + n^{2} \left(32112 \sigma^{8} - 6948 \sigma^{6} + 246 \sigma^{4}\right) + n \left(- 20160 \sigma^{8} + 5040 \sigma^{6} - 252 \sigma^{4} + \sigma^{2}\right)\right)$ \\
10 &  $945 n^{5} \sigma^{10} + n^{4} \left(- 44100 \sigma^{10} + 9450 \sigma^{8}\right) + n^{3} \left(303660 \sigma^{10} - 99120 \sigma^{8} + 6825 \sigma^{6}\right) + n^{2} \left(- 623376 \sigma^{10} + 240840 \sigma^{8} - 24438 \sigma^{6} + 501 \sigma^{4}\right) + n \left(362880 \sigma^{10} - 151200 \sigma^{8} + 17640 \sigma^{6} - 510 \sigma^{4} + \sigma^{2}\right)$ \\
\bottomrule
\end{tabular}

}
\caption{Central Moments of Binomial Distribution. Above we denote $\sigma^2=p(1-p)$.}
\label{tab:moments}
\end{table}

\begin{figure}
\centering
\begin{tikzpicture}
    \begin{axis}[view={-45}{45},xlabel={$n$},ylabel={$p$},zlabel={$\mathbb{E}[S-\mathbb{E}[S]]^d$}]
    \addplot3[surf,colormap = {whiteblack}{color(0cm)  = (white);color(0.5cm) = (black!75)},domain=0:99] table[col sep=comma] {binomial_moments.txt};
    \end{axis}
\end{tikzpicture}
\caption{The $d$-th central moment of $S(n,p)$, where $d=6$.}
\end{figure}
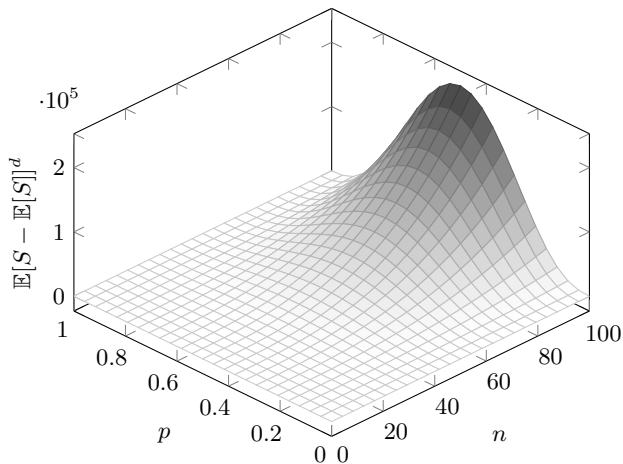

\paragraph{Positive Polynomial Representation}

As mentioned in the introduction, the only closed-form formula due to~\cite{knoblauch2008closed} is an alternating
sum with no readable leading term, which makes it hard to use; in particular the discussion in~\cite{knoblauch2008closed} fails to give non-trivial bounds on binomial moments. The novelty of our work is a formula consisting of \emph{positive terms}. This makes it more stable for numerical computations and more handy in theoretical analysis.
\begin{theorem}[Stable Expression]\label{cor:binom_central}
For $S\sim\mathrm{Binom}(n,p)$, $q\triangleq 1-p$ and any positive integer $d$ the following holds
\begin{align*}
\mathbb{E}[(S-\mathbb{E}[S])^d] = \sum_{k=1}^{\lfloor\frac{d}{2}\rfloor}\binom{n}{k} (pq)^{k}{\sum_{d_1\ldots d_k\geqslant 2} \binom{d}{d_1\ldots d_k} \prod_{i=1}^{k} (q^{d_i-1}-(-p)^{d_i-1})}.
\end{align*}
\end{theorem}

\begin{remark}[Non-negativity]
The terms under the sum are all non-negative when $p\leqslant \frac{1}{2}$ and are all negative when $p>\frac{1}{2}$.
Since $n-S\sim\mathrm{Binom}(n,q)$, it follows that $S-\mathbb{E}[S]\sim -(\mathrm{Binom}(n,q)-qn)$; with the help of this identity, studying the central binomial moments can be always reduced to the case $p<\frac{1}{2}$.
\end{remark}

\begin{remark}[Closed-form Symmetric Formula]
The above result gives an explicit formula for \Cref{thm:structure}, and provides an alternative proof of that result.
\end{remark}

\subsubsection{Asymptoticaly Sharp Moment Bounds}

To illustrate how useful is our \emph{positive} representation established in \Cref{cor:binom_central}, we derive the \emph{sharp bounds} on (normalized) central binomial moments. This problem has remained open so far; some recent works used adhoc upper bounds to estimate the binomial moments and tails (works on random projections, particularly~\cite{jagadeesan2019understanding}).

\begin{theorem}\label{thm:asymptotic}
Let $S\sim\mathrm{Binom}(n,p)$. Then for any positive even $d$ we have
\begin{align*}
\mathbb{E}[(S-\mathbb{E}[S])^d]^{1/d}=\Theta(1)\cdot\max\left\{k^{1-\frac{k}{d}}\cdot (n\sigma^2)^{\frac{k}{d}}:   k = 2\ldots d/2\right\}.
\end{align*}
\end{theorem}
The estimate is uniformly sharp in all parameters; for the special case when $n\sigma^2\to \infty $ and $d$ is fixed, the $d$-th central moment grows as $\Theta(n\sigma^2)^{d/2}$ which matches the central limit theorem combined with the explicit formulas for moments of the normal distribution~\cite{winkelbauer2012moments}. In contrast to \Cref{thm:asymptotic}, the formula in ~\cite{knoblauch2008closed} gives in this setup only much worse $O(nq)^{d}$, which anyway follows trivially since $S$ is bounded by $n$. The proof uses \Cref{cor:binom_central} and Stirling's formulas to estimate the growth of leading terms (see \Cref{proof:asymptotic}).

\subsubsection{Algorithms and Implementation}

We have seen that the variance-based representation in \Cref{thm:structure} is useful,
yet it is not immediate \emph{how to compute} this polynomial.
To this end, we develop two different algorithms, both implemented in the popular Python symbolic algebra package \texttt{Sympy}~\cite{sympy}.

 \Cref{alg:1} takes advantage of the \emph{fumdamental theorem on symmetric polynomials} (see for example \cite{gaudry2006evaluation}).
Specifically, there is an explicit procedure for converting any symmetric polynomial in $p,q$ into a polynomial 
in variables $p+q,pq$ (the basic symmetric polynomials in two variables); substituting $p+q=1$ we are left with a polynomial in $\sigma^2=pq$. For even $d$ we start with a symmetric polynomial and obtain a polynomial in $\sigma^2$. In turn,
for odd $d$ we start with an anti-symmetric polynomial and apply this procedure to its symmetric part, which results in a polynomial in $\sigma^2$ plus the factor $q-p=1-2p$. 
\begin{algorithm}[h!]
  \scriptsize
\SetAlgoLined
\KwResult{The $d$-th central binomial moment in terms of $n$ and $\sigma^2=p(1-p)$. \\ The outcome is a polynomial with integer coefficients.}
 $U\gets \mathbb{E}(S-\mathbb{E}[S])^d$, with $U\in\mathbb{Z}[n][p,q]$ \tcp{write as polynomial in $p,q$, e.g. by \Cref{cor:binom_central}}
\If{$d$ even}{
$V\gets \texttt{symmetrize}(U),\ V\in \mathbb{Z}[n][pq,p+q]$ \tcp{represent by elementary symmetric polynomials}
$V\gets V|_{p+q\gets 1}$ \tcp{substitute the relation $p+q=1$}
}

\Else{
$V\gets V/(q-p)$  \tcp{subtract 'unsymmetric' part}
$V\gets \texttt{symmetrize}(U),\ V\in \mathbb{Z}[n][pq,p+q]$ \tcp{represent by elementary symmetric polynomials}
$V\gets V|_{p+q\gets 1}$ \tcp{substitute the relation $p+q=1$}
$V\gets (1-2p)\cdot V$ \tcp{add unsymmetric part back}
}
\Return{V}
 \caption{Variance Formula for Central Binomial Moments}
\label{alg:1}
\end{algorithm}

In turn \Cref{alg:2} uses the power of \emph{elimination theory}, to recover the representation 
in $\sigma^2$ from any formula written in terms of $p$. Essentially, it simplifies the input polynomial in $p$ with respect to the polynomial $\sigma^2-p(1-p)$
leveraging the elimination properties of \emph{Groebner bases} (see for example~\cite{buchberger2001grobner}). The output is a polynomial in $\sigma^2$ (plus the factor $1-2p$ for odd $d$).
\begin{algorithm}[h!]
  \scriptsize
\SetAlgoLined
\KwResult{The $d$-th central binomial moment in terms of $n$ and $\sigma^2=p(1-p)$.  \\
The output is a polynomial with integer coefficients.}
 $U\gets \mathbb{E}(S-\mathbb{E}[S])^d,\ U\in\mathbb{Z}[n,p]$ \tcp{represent as polynomial in $n,p$, for example  
\Cref{cor:binom_central} and substitution $q\gets 1-p$}
\If{$d$ even}{
$F_1,F_2\gets U,\sigma^2-p(1-p)$ \tcp{moment and variance formulas}
$G \gets \texttt{GroebnerB}(\textrm{poly}=\{F_1,F_2\},\textrm{vars}=(p,n,\sigma^2),\textrm{order}=\textrm{lex})$ \tcp{Groebner basis, lex order}
$V\gets G\cap \mathbb{Z}[n,\sigma^2]$ \tcp{extract variance-dependent formula}
}
\Else{
 $F_1,F_2\gets U/(1-2p),\sigma^2-p(1-p)$ \tcp{subtract 'odd' part}
$G \gets \texttt{GroebnerB}(\textrm{poly}=\{F_1,F_2\},\textrm{vars}=(p,n,\sigma^2),\textrm{order}=\textrm{lex})$ \tcp{Groebner basis, lex order}
$V\gets G\cap \mathbb{Z}[n,\sigma^2]$ \tcp{extract variance-dependent part}
$V\gets (1-2p)\cdot V$ \tcp{add 'odd' part back}
}
\Return{V}
 \caption{Variance Formula for Central Binomial Moments}
\label{alg:2}
\end{algorithm}

\subsection{Organization}
The remainder of the paper is organized as follows:
\Cref{sec:proofs} gives proofs of the results, \Cref{sec:app} discusses the application to sharp asymptotics, 
\Cref{sec:implementation} presents the Python implementation and finally \Cref{sec:conclusion} concludes the work.

\section{Proofs}\label{sec:proofs}

\subsection{First Proof of \Cref{cor:binom_raw}}\label{proof:binom_raw:factorial}

The proof is based on the fact that the \emph{factorial moments} of the binomial distribution are easy to compute. Namely (see \cite{potts1953note,ewis2012central}) we have
\begin{proposition}[Factorial Moments of Binomial Distribution]\label{prop:factorial_binom_moments}
Let $S\sim \mathrm{Binom}(n,p)$. Then the following holds
\begin{align}\label{eq:binomial_fact_moment}
\mathbb{E}[S^{\underline{d}}] = \binom{n}{d} p^d.
\end{align}
\end{proposition}

Then it remains to connect factorial moments to standard moments, or in other terms: factorial powers to powers. It is well-known (see for example the discussion in~\cite{boyadzhiev2012close}) that this base change is given in terms of the Stirling numbers of the second kind. We state this fact formally below
\begin{proposition}[Base Change from Falling Powers to Powers]\label{prop:stirling_base_change}
For positive integers $x$ and $d$ the following holds
\begin{align}\label{eq:stirling_change_base}
x^d = \sum_{k=0}^{d}\stirling{d}{k} x^{\underline{k}}.
\end{align}
\end{proposition}
Now \Cref{cor:binom_raw} follows by applying \Cref{prop:stirling_base_change} to $x:= S$, and then using \Cref{prop:factorial_binom_moments} to compute the expectation of $S^{\underline{k}}$ for $k=0,\ldots,d$.

\subsection{Second Proof of \Cref{cor:binom_raw}}\label{proof:binom_raw:direct}

Here we take a direct approach, writing $S = \sum_{i=1}^{n}X_i$ where $X_i \sim^{iid}\mathrm{Bern}(p)$.
Using the multinomial expansion and the independence of $X_i$ we obtain
\begin{align*}
\mathbb{E}[S^d] = \sum_{d_1,\ldots,d_n} \binom{d}{d_1\ldots d_n}\prod_{i}\mathbb{E}[X_i^{d_i}].
\end{align*}
We now group the expressions in the above sum, depending on the number of non-zero elements in $(d_i)_i$. Denoting
$\|(d_i)\|_0 = \#\{i: d_i > 0\}$ and using the fact that $\mathbb{E}X_i^{d_i} = p$ when $d_i>0$ we obtain
\begin{align*}
\mathbb{E}[S^d] = \sum_{k=0}^{d}\sum_{(d_i)_i: \|(d_i)_i\|_0 = k} \binom{d}{d_1\ldots d_n}p^k.
\end{align*}
By the symmetry of multinomial coefficients this equals
\begin{align*}
\mathbb{E}[S^d] = \sum_{k=0}^{d}\sum_{d_1\ldots d_k>0} \binom{d}{d_1\ldots d_k} \binom{n}{k} p^k.
\end{align*}
Finally, we observe that the expression $\sum_{d_1\ldots d_k>0} \binom{d}{d_1\ldots d_k}$
counts the number of ways of partitioning $\{1,\ldots,d\}$ into $k$ non-empty \emph{labeled} subsets;
thus, this numbers equals $k!\cdot\stirling{d}{k}$ which finishes the proof.

\subsection{First Proof of \Cref{thm:structure}}\label{proof:thm:structure:algebra}

From \Cref{eq:density} we obtain
\begin{align*}
\mathbb{E}[(S-\mathbb{E}[S])^d] = \sum_{k} \binom{n}{k}p^k q^{n-k} (k-np)^d.
\end{align*}
Replacing $k$ by $n-k$ and using the symmetry of binomial coefficients we obtain
\begin{align*}
\mathbb{E}[(S-\mathbb{E}[S])^d] = \sum_{k} \binom{n}{k}p^{n-k} q^{k} (nq-k)^d.
\end{align*}
When $d$ is even, comparing these two equivalent expressions we see that they are symmetric as polynomials in $p$ and $q$.
By the fundamental theorem of symmetric polynomials, this can be written as a polynomial in $pq$ and $p+q$; in our case $p+q=1$ and the claim follows. If $d$ is odd then $(nq-k)^d = -(k-nq)^d$ and we get anti-symmetric polynomials in $p,q$ which can be written as a product of $p-q$ and a symmetric polynomial. The latter, by the fundamental theorem, is a polynomial in $p+q$ and $pq$; since $p+q=1$ the result follows.

\subsection{Second Proof of \Cref{thm:structure}}\label{proof:thm:structure:combinatorics}
By inspecting the products $\prod_{i=1}^{k}(q^{d_i-1}-(-p)^{d_i-1}$ that appear in
 \Cref{cor:binom_central} it can be seen that each of them is symmetric in $p,q$ when $\sum_i d_i = d$ is even, and anti-symmetric when
$\sum_i d_i = d$ is odd. This is because $p,q\to q^{d_i-1}-(-p)^{d_i-1}$ is symmetric when $d_i$ is even and anti-symmetric otherwise.
The claim now follows.

\subsection{Proof of \Cref{cor:binom_central}}

As in the proof of \Cref{cor:binom_raw} we arrive at
\begin{align*}
\mathbb{E}[(S-\mathbb{E}[S])^d] = \sum_{k=0}^{d}\sum_{d_1\ldots d_k>0} \binom{n}{k} \binom{d}{d_1\ldots d_k}\prod_{i=1}^{k}
\mathbb{E}[(X_i-\mathbb{E}[X_i])^{d_i}].
\end{align*}
Denote $x = 1-\frac{1}{1-p}$, then
$
\mathbb{E}[(X_i-\mathbb{E}[X_i])]^{d_i} = p(1-p)^{d_i}(1-x^{d_i-1})
$ and thus
\begin{align*}
\mathbb{E}[(S-\mathbb{E}[S])^d] = (1-p)^d\sum_{k=0}^{d}\binom{n}{k} p^{k} \sum_{d_1\ldots d_k>0} \binom{d}{d_1\ldots d_k} \prod_{i=1}^{k} (1-x^{d_i-1}).
\end{align*}
With some further simplifications and grouping we can write
\begin{align*}
\mathbb{E}[(S-\mathbb{E}[S])^d] = (1-p)^d\sum_{k=1}^{\lfloor\frac{d}{2}\rfloor}\binom{n}{k} p^{k}\underbrace{\sum_{d_1\ldots d_k\geqslant 2} \binom{d}{d_1\ldots d_k} \prod_{i=1}^{k} (1-x^{d_i-1})}_{U_k},
\end{align*}
or equivalently
\begin{align*}
\mathbb{E}[(S-\mathbb{E}[S])^d] = \sum_{k=1}^{\lfloor\frac{d}{2}\rfloor}\binom{n}{k} (pq)^{k}{\sum_{d_1\ldots d_k>1} \binom{d}{d_1\ldots d_k} \prod_{i=1}^{k} (q^{d_i-1}-(-p)^{d_i-1})}.
\end{align*}
This finishes the proof. In addition to that, in what follows, we discuss how to further group terms and speed up computations.
We can write
\begin{align*}
U_k=\sum_{j=0}^{k}\binom{k}{j} \sum_{d_1\ldots d_{j}>1}\sum_{d_{j+1}\ldots d_k>1} \binom{d}{d_1\ldots d_k}(-1)^{j}x^{d_1+\ldots+d_{j}-j}
\end{align*}
since $\binom{d}{d_1\ldots d_k}=\binom{d_1+\ldots+d_j}{d_1\ldots d_j}\cdot \binom{d-(d_1+\ldots+d_j)}{d_{j+1}\ldots d_k}\cdot \binom{d}{d_1+\ldots+d_j}$ we obtain
\begin{align*}
U_k = \sum_{j=0}^{k}\binom{k}{j}\sum_{\ell}\binom{d}{\ell}\sum_{d_1\ldots d_{k}>1}\binom{\ell}{d_1\ldots d_{j}}
\binom{d-\ell}{d_{j+1}\ldots d_{k}}(-1)^{j} x^{\ell-j},
\end{align*}
and thus
\begin{align*}
U_k &= \sum_{j=0}^{k}\sum_{\ell=0}^{d}\binom{k}{j}\binom{d}{\ell} j! (k-j)! S_2(\ell,j)S_2(d-\ell,k-j) (-1)^{j} x^{\ell-j} \\
& = k!\sum_{\ell=0}^{d}\binom{d}{\ell}\sum_{j=0}^{k} S_2(\ell,j)S_2(d-\ell,k-j) (-1)^{j} x^{\ell-j},
\end{align*}
where $S_2(n,k)$ denotes the number of ways of partitioning an $n$-element set into $k$ subsets of cardinality at least $2$ (a variation on Stirling numbers of the second kind). This can be used to develop an equivalent, but faster to compute, formula.

\section{Applications}\label{sec:app}

\subsection{Proof of \Cref{thm:asymptotic}}\label{proof:asymptotic}

When $p\leqslant \frac{1}{2}$ we have $q\geqslant p$ and thus
$0\leqslant q^{d_i-1}-(-p)^{d_i-1}\leqslant p+q=1$  for any $d_i\geqslant 2$. In view of \Cref{cor:binom_central}, we obtain the following bound
\begin{align*}
\mathbb{E}[(S-\mathbb{E}[S])^d] \leqslant \sum_{k=1}^{\lfloor\frac{d}{2}\rfloor}\binom{n}{k} (pq)^{k}{\sum_{d_1\ldots d_k\geqslant 2} \binom{d}{d_1\ldots d_k}}.
\end{align*}
Since we have
\begin{align*}
\sum_{d_1\ldots d_k\geqslant 2} \binom{d}{d_1\ldots d_k}
\leqslant \sum_{d_1\ldots d_k\geqslant 0} \binom{d}{d_1\ldots d_k} = k^d
\end{align*}
we further obtain
\begin{align*}
\mathbb{E}[(S-\mathbb{E}[S])^d] \leqslant \sum_{k=1}^{\lfloor\frac{d}{2}\rfloor}\binom{n}{k} (pq)^{k}k^d.
\end{align*}
Denoting $\sigma^2=pq$, using the elementary bound $\binom{n}{k}\leqslant (n/k)^k$, setting $r=d/k$ and using the asymptotic $r^{1/r}=\Theta( 1)$ we finally obtain
\begin{align}
    \mathbb{E}[(S-\mathbb{E}[S])^d] \leqslant O(1)\cdot \max\{k^{1-\frac{k}{d}}\cdot (n\sigma^2)^{\frac{k}{d}}:  2\leqslant k \leqslant \lfloor d/2\rfloor\}.
\end{align}
We now move on to the lower bound. When $d_i$ are even, we have $q^{d_i-1}-(-p)^{d_i-1} = q^{d_i-1}+p^{d_i-1}\geqslant 2\cdot \frac{1}{2^{d_i-1}}$ by Jensen's inequality applied to the function $u\to u^{d_i-1}$ and $p+q=1$. Since in the summation we consider $(d_i)_i$ such that $\sum_{i}d_i=d$ and $d_i\geqslant 2$, we obtain
\begin{align*}
    \prod_{i=1}^k(q^{d_i-1}-(-p)^{d_i-1}) \geqslant \prod_{k}2^{2-d_i} = 2^{2k-d}.
\end{align*}
Let us write $d = k\cdot r + \ell$ with non-negative even integers $\ell,r$ such that $k\leqslant \ell<2k$; this is possible when $d$ is even, 
by dividing with the remainder $d = k\cdot r + \ell$ where $0\leqslant \ell <k$ and replacing $r:=r-1,\ell:=\ell+k$ when $r$ is odd. 
Define $d_i = r+2$ when $i\leqslant \ell/2$ and $d_i = r$ when $\ell/2<i\leqslant k$.
Using Stirling's approximation $m! = \Theta(m)^{m+\frac{1}{2}}$, and $r = \Theta(d/k)$ we obtain
\begin{align*}
\binom{d}{d_1\ldots d_k} =\frac{\Theta(d)^{d+1/2}}{ \Theta(r)^{\sum_{i=1}^k d_i +k/2} }= \frac{\Theta(d)^{d+1/2}}{\Theta(r)^{d+k/2}}
= \Theta(d/r)^{d} = \Theta(k)^{d}.
\end{align*}
The above two bounds, in view of \Cref{cor:binom_central}, imply that
\begin{align*}
\mathbb{E}[(S-\mathbb{E}[S])^d] \geqslant \sum_{k=1}^{\lfloor\frac{d}{2}\rfloor}\binom{n}{k} (pq)^{k}  \Omega(k)^d.
\end{align*}
Denoting $\sigma^2=pq$, using the elementary bound $\binom{n}{k}\geqslant (n/\mathrm{e}k)^k$, setting $r=d/k$ and using the asymptotic $r^{1/r}=\Theta( 1)$ we finally obtain
\begin{align}
    \mathbb{E}[(S-\mathbb{E}[S])^d] \geqslant \Omega(1)\cdot \max\{k^{1-\frac{k}{d}}\cdot (n\sigma^2)^{\frac{k}{d}}:  2\leqslant k \leqslant \lfloor d/2\rfloor\}
\end{align}
which finishes the proof.

\section{Conclusion}\label{sec:conclusion}

This paper introduces novel and simpler formulas for binomial moments, derived by a combinatorial argument coupled with clever algebraic simplification which relies on symmetrization. We show applications of independent interest, such as deriving sharp asymptotics for the growth of central binomial moments. Moreover, explicit algorithms and the working implementation are provided.

%
%
%
\bibliographystyle{apalike}
\bibliography{citations}

\appendix

\section{Implementation}\label{sec:implementation}

\definecolor{backcolour}{rgb}{0.95,0.95,0.92}
\definecolor{codegreen}{rgb}{0,0.6,0}
\definecolor{mymauve}{rgb}{0.58,0,0.82}
\lstdefinestyle{myStyle}{
    backgroundcolor=\color{backcolour},   
    keywordstyle=\color{blue},
    commentstyle=\color{codegreen},
    stringstyle=\color{mymauve},
    basicstyle=\footnotesize,
    breakatwhitespace=false,         
    breaklines=true,                 
    keepspaces=true,                 
    xleftmargin=10pt,
    showspaces=false,                
    showstringspaces=false,
    showtabs=false,                  
    tabsize=4,
    frame=L,
    captionpos=b, 
    abovecaptionskip=1\baselineskip,
}
\lstset{style=myStyle}

\lstset{language=Python,basicstyle=\scriptsize}          

\begin{lstlisting}[frame=single,caption=Preliminaries,float]
import itertools
import sympy as sm
from sympy import Symbol
from sympy import polys
from sympy.functions.combinatorial.numbers import stirling,binomial
from sympy.functions.combinatorial.numbers import factorial
from sympy.functions.combinatorial.factorials import FallingFactorial

def multinomial_coef(n,ks):
  if n!=sum(ks):
    return 0
  elif len(ks)==1:
    return 1
  else:
    return binomial(n,ks[0])*multinomial_coef(n-ks[0],ks[1:])
\end{lstlisting}

\begin{lstlisting}[frame=single,caption={Stable Formula for Central Moments},float]
def central_binom_moment(d=2):
  ''' output as poly in trials number n, success prob. p and q=1-p '''
  n = Symbol('n')
  p = Symbol('p')
  q = Symbol('q')
  out = 0
  for k in range(1,int(d/2)+1):
    tmp = 0
    for dks in itertools.product(range(2,d+1),repeat=k):
      polyx = multinomial_coef(d,dks) 
      for dk in dks:
        polyx = polyx * (q**(dk-1)-(-p)**(dk-1))
      tmp = tmp + polyx
    out = out + binomial(n,k) * (p*q)**k * tmp

  return out.simplify()
\end{lstlisting}

\begin{lstlisting}[frame=single,caption={Variance Formula for Central Moments by Variable Elimination},float]
def central_binom_moment_pretty1(d=2):
  ''' output as polynomial in number of trials and success variance '''

  out = central_binom_moment(d=d).combsimp()
  p = Symbol('p')
  q = Symbol('q')
  s = Symbol('sigma')
  n = Symbol('n')
  out = out.subs(q,1-p)
  if d % 2 == 0:
    out = polys.groebner([out,s**2-p*(1-p)],p,n,s,order='lex')[-1]
  elif d % 2 == 1:
    out = polys.div(out,1-2*p)[0]
    out = polys.groebner([out,s**2-p*(1-p)],p,n,s,order='lex')[-1]
    out = (1-2*p)*out
  return out
\end{lstlisting}

\begin{lstlisting}[frame=single,caption={Variance Formula by Symmetrization},float]
def central_binom_moment_pretty2(d=2):
  ''' output as polynomial in number of trials and success variance '''

  out = central_binom_moment(d=d)
  p,q = Symbol('p'),Symbol('q')
  s = Symbol('sigma')
  if d % 2 == 0:
    out = polys.symmetrize(out,gens=[p,q])[0]
  elif d % 2 == 1:
    out = polys.div(out,1-2*p)[0]
    out = polys.symmetrize(out,gens=[p,q])[0]
    out = (1-2*p)*out
  out = out.subs(p+q,1)
  out = out.subs(p*q,s**2)
  return out
\end{lstlisting}

\end{document}